\documentclass{AIMS}
 \usepackage{amsmath}
 \usepackage{epsfig} 
 \usepackage{graphics} 

\newtheorem{theorem}{Theorem}

\theoremstyle{definition}

\title[How distant is the ideal filter of being a causal one?]
      {How distant is the ideal filter of being a causal one?}

\author[J. M. Almira and A. E. Romero ]{}


\email{jmalmira@ujaen.es; aeromero@decsai.ugr.es}

\begin{document}

\maketitle

\centerline{\scshape J. M. Almira}
\medskip
{\footnotesize{}
 \centerline{Departmento de Matem\'{a}ticas}
  \centerline{Universidad de Ja\'{e}n, E.P.S. Linares}
   \centerline{23700 Linares (Ja\'{e}n) Spain}
} 

\medskip

\centerline{\scshape A. E. Romero}
\medskip
{\footnotesize{}
 \centerline{Departamento de Ciencias de la Computaci\'on e Inteligencia Artificial}
  \centerline{Universidad de Granada. Facultad de Ciencias. Campus Fuentenueva}
   \centerline{18071 Granada.  Spain}
} 

\medskip

  \medskip
\begin{abstract}
In this paper the characterization as convolution operators of filters sending finite energy signals to bounded signals   is used to prove several theoretical results concerning the distance between the ideal filter and the spaces of physically realizable filters. Both the analog and the digital cases are studied and the formulas for the distance and the angle between the filters in each case are also given.
\end{abstract}

\section{Introduction and motivation}
One of the main principles on which the deterministic mathematical physics is based is the statement that Nature obeys certain universal laws, and these laws can be usually expressed in terms of partial differential equations. Moreover, this view is commonly extended assuming some additional properties. For example, there is a strong feeling supporting that Nature is uniform. This means that a law that holds true at a certain place and a certain moment holds true everywhere and all times. In particular, this implies that the involved partial differential equations in the description of the physical laws should have coefficients independent on the time. Moreover, under certain additional conditions, many of these laws are linear in nature. Hence, the continuous linear operators defined between two function spaces that have the property of being translation-invariant are of great interest in Applied Mathematics. These operators are usually named filters by the engineering community (See
\cite[p. 14 and Definition 34.1.1]{witomsky}).

More precisely, let $X$ and $Y$ be two function spaces with the property that their elements are functions of the real variables $t_1,\dots,t_n$ and, for each $h\in\mathbb{R}^n$, let $\tau_{h}:\mathbb{R}^n\longrightarrow \mathbb{R}^n$ be the translation operator $\tau_{h}(a)=a-h$. By definition, an analog filter with input signals from $X$ and output signals from $Y$  is a continuous linear operator $L:X\longrightarrow Y$ such that, for every  $h\in\mathbb{R}^n$,  the relation $L\left(x(\tau_h(s))\right)(t)=L(x(s))(\tau_h(t))$ holds true. On the other hand, if $X,Y$ are sequence spaces (i.e. their elements are $n$-dimensional sequences $\{a_k\}_{k\in\mathbb{Z}^n}$), we say that  $L:X\longrightarrow Y$ is a digital filter if it is a continuous linear operator and the relation $L\left(\{x_k\}_{k\in\mathbb{Z}^n}\right)=\{y_k\}_{k\in\mathbb{Z}^n}$ implies the relation $L\left(\{x_{k-N}\}_{k\in\mathbb{Z}^n}\right)=\{y_{k-N}\}_{k\in\mathbb{Z}^n}$ for every $N\in\mathbb{Z}^n$. The best known examples of filters are the convolution operators, which are operators of the form
\[
L_h(x)(t)=(x\ast h)(t)=\int_{-\infty}^{\; \infty}x(s)h(t-s)ds.
\]
Indeed, the name ``filter'' has its origin in the fact that under quite general conditions the Fourier transform of a convolution is a product, so that the filters $L_h$ can be represented in the frequency domain as $Y=X\cdot H$, where $X=\mathcal{F}(x)$, $H=\mathcal{F}(h)$ and $Y=\mathcal{F}(L(x))$. Hence, if the function $H(\xi)$ satisfies $H(\xi)=1$ for the frequencies $\xi\in W_1$ and $H(\xi)=0$ for $\xi\in W_2$ then we understand that $L_h$ allows the frequency content of the signal $x$ that belongs to the set $W_1\subseteq \mathbb{R}$ but removes the frequency content belonging to $W_2\subseteq \mathbb{R}$, so that we may interpret that $L_h$ filters some frequency contents of $x$.

Filter Theory is an interesting branch of Mathematical Analysis that contains many beautiful results. For example, the characterization of filters between some kinds of function spaces $X$, $Y$ strongly depends on the geometric and/or analytical properties of these spaces, and it is usually a difficult problem (See \cite{Hormander} for a classical paper, \cite[Chapter zero]{hytonenthesis} for an historical overview and \cite{hytonen1}, \cite{hytonen2} for some recent results. For distributions, a classical result by L. Schwartz stating that every filter $T:\mathcal{E}'\to\mathcal{D}'$ is a convolution operator (where $\mathcal{E}'$ denotes the space of distributions with compact support and $\mathcal{D}'$ is the space of distributions), is to be found in \cite[Theorem 5.8.1, page 144]{zemanian}). Moreover, this theory is clearly related to such diverse subjects as Fourier Analysis, Functional Analysis and Signal Processing. From now on, we will restrict our attention to the one-dimensional case ($n=1$) because it is precisely in that context where physically realizable and causal filters have are meaningful. Under this assumption, it is quite natural to identify our unique independent variable, $t$,  with time.

Let us accord that a filter $L$ is physically realizable if there exists a certain \linebreak $T\geq 0$ such that for every time value $t\in\mathbb{R}$ the output of the filter at $t$ depends on the values of the input signal on the interval $(-\infty,t+T]$. This means that (at least theoretically) the filter can be realized under the understanding that we allow some delay time $T$. The filter does not depend on all the past and future values of the signal but just on the past and perhaps some future values limited in time. When $T=0$ we say that the filter is causal. The concept of causal filter was introduced in the literature by N. Wiener in 1926 in a paper devoted to give an adequate interpretation of Heaviside's operational calculus \cite{Wienerheaviside}. That paper is also famous because it was an inspiration point for the introduction by L. Schwartz of his theory of distributions.

Now, it is well known that the analog ideal filter, which is given in the frequency domain by
\begin{equation}\label{primera}
L(X)(\xi)=X(\xi)H(\xi),
\end{equation}
where
$H=\chi_{_{[a,b]}}$  and $\chi_{_{[a,b]}}(\xi)=1$ for $\xi\in[a,b]$,  $\chi_{_{[a,b]}}(\xi)=0$ for $\xi \not \in [a,b]$, is not physically realizable. The simplest reason is that the inverse Fourier transform of $\chi_{_{[a,b]}}$ does not vanish over any interval of the form $(-\infty,-T)$.  Another explanation appears when we take into account an important theorem by Paley and Wiener \cite[Theorem XII, p. 16]{PW} \cite[p. 35]{Paarmann} where they characterize physically realizable filters in the frequency domain as those of the form $(\ref{primera})$ for which $H(\xi)\in L^2(\mathbb{R})$ and
\[
\int_{-\infty}^{\; \infty} \frac{|\log |H(\xi)|\, |}{1+\xi^2}\, d\xi<\infty.
\]
(With such a characterization it is obvious that $H(\xi)$ cannot vanish on any open set of the real line!). Wiener was so proud of this characterization that he commented the result several times in printed form. For example, in \cite[p. 37]{Wieneryellowperil} he said:
\medskip

{\footnotesize {\it This [result] plays a very important part in the theory of filters. It states that, in any realizable network whatever, the attenuation, taken as a function of the frequency $\omega$, and divided by $1+\omega^2$, yields an absolutely integrable function of the frequency. This results from the fact that the attenuation is the logarithm of the absolute value of the transform of [the response to the unit impulse] $f(t)$ which vanishes for negative $t$; or, in other words, because strictly no network can foretell the future. Thus no filter can have infinite attenuation in any finite band. The perfect filter is physically unrealizable by its very nature, not merely because of the paucity of means at our disposal. No instrument acting solely on the past has a sufficiently sharp discrimination to separate one frequency from another with unfailing accuracy.}}
\medskip

Moreover, in his autobiography \cite[p. 168]{Wienerautobiography}, when speaking about his mathematical work with Paley, he said:
\medskip

{\footnotesize {\it One interesting problem which we attacked together was that of the conditions restricting the Fourier transform of a function vanishing on the half line. This is a sound mathematical problem on its own merits, and Paley attacked it with vigor, but what helped me and did not help Paley was that it is essentially a problem in electrical engineering. It had been known for many years that there is a certain limitation on the sharpness with which an electric wave filter cuts a frequency band off, but the physicists and engineers had been quite unaware of the deep mathematical grounds for these limitation. In solving what  was for Paley a beautiful and difficult chess problem, completely contained within itself, I showed at the same time that the limitations under which the electrical engineers were working were precisely those which prevent the future from influencing the past.}}
\medskip

The result was, moreover, a key step for the proof of several fundamental theorems in Harmonic Analysis such as Carleman's characterization of quasi-analytic functions. In our opinion, Wiener was in his own right to be proud of his result and all the ``philosophical''  interpretations he gives to it are essentially correct and illuminating. We have wondered if a certain quantitative estimation of the ``far away'' ideal filters are of being physically realizable already exists. We have not found such an estimation in the literature. The main goal of this paper is to make several computations in this direction. In order to give a precise focus to our computations we first introduce a result which characterizes certain filters as convolution filters. After that, we dedicate a section to the analog case and another to the digital one.

\section{Characterization of filters as convolution operators}

Let $\mathbf{B_a}=\mathbf{B}(L^2(\mathbb{R}),$ $L^{\infty}(\mathbb{R}))$ be the normed space of linear bounded operators
$L:L^2(\mathbb{R})\longrightarrow L^{\infty}(\mathbb{R})$, and let us consider $\mathbf{F_a}\subset \mathbf{B_a}$ be the subspace of analog filters. Thus, $L$ belongs to $\mathbf{F_a}$ if it sends finite energy analog signals to bounded analog signals, is a bounded linear operator and  is time invariant. Analogously, we can define $\mathbf{B_d}=\mathbf{B}(\ell^2(\mathbb{Z}),\ell^{\infty}(\mathbb{Z}))$, the normed space of linear bounded operators $L:\ell^2(\mathbb{Z})\longrightarrow \ell^{\infty}(\mathbb{Z})$, and $\mathbf{F_d}\subset \mathbf{B_d}$, the subspace of digital filters. Thus, $L$ belongs to $\mathbf{F_d}$ if it sends finite energy digital signals to bounded digital signals, is a bounded linear operator and is time invariant. We devote this section to characterize the elements of $\mathbf{F_a}$ and $\mathbf{F_d}$ as convolution operators in such a way that both spaces are naturally doted of a Hilbert space structure. These results are not new (see, for example, \cite{hewittross}, \cite{Unni}, where several generalizations of them are proved). By the contrary, they are well know but they are also not attributable to any person. Moreover, they are not available in the literature in the simple version we have stated here but they are just special cases of quite difficult results. Thus we think a direct proof will make them more visible and useful.

\begin{theorem}The map $\phi:L^2(\mathbb{R})\longrightarrow \mathbf{F_a}$ given by $\phi(h)=L_h$, where $L_h(x)=x\ast h$, is an isometry. In particular, $\mathbf{F_a}$, with its usual norm (inherited from $\mathbf{B_a}$) is a Hilbert space with inner product given by $\langle L_h,L_g \rangle = \langle h,g \rangle$.
\end{theorem}

\noindent \textbf{Proof.} Obviously, $\phi$ is a linear injective map but, in principle, there is the possibility that some elements $L$ of $\mathbf{F_a}$ are not of the form $L_h$, for a certain $h\in L^2(\mathbb{R})$. Thus, we must prove that $\phi$ is surjective and preserves the norm.
 Given $L\in\mathbf{F_a}$ and $t\in\mathbb{R}$, the map $L_t:L^2(\mathbb{R})\longrightarrow\mathbb{C}$, defined by $L_t(x)=L(x)(t)$, is a linear bounded functional. It follows from the Riesz representation theorem \cite[p. 24]{krantz} that there is $y_t\in L^2(\mathbb{R})$ such that $L_t(x)=\langle x,y_t \rangle$, for any $x\in L^2(\mathbb{R})$. Let us denote $y(t,s)=y_t(s)$. Time invariance of $L$ means that for every $K\in\mathbb{R}$ and every $x\in L^2(\mathbb{R})$,
\[
L\left(x(\cdot+K)\right)(t)=L(x)(t+K).
\]
Hence
\[
\langle x(\cdot+K),y_t\rangle = \langle x,y_{t+K} \rangle,
\]
so that
\[
\int_{-\infty}^{\;\infty}x(s+K)\,\overline{y(t,s)}\,ds=\int_{-\infty}^{\;\infty}x(\tau)\,\overline{y(t,\tau-K)}\,d\tau=\int_{-\infty}^{\;\infty}x(s)\,\overline{y(t+K,s)}\,ds.
\]
This implies that $y(t,s-K)=y(t+K,s)$. If we take $h(t)=\overline{y(t,0)}$, then
\[
y(t,s)=y(t,0-(-s))=y(t-s,0)=\overline{h(t-s)}.
\]
Hence $L=L_{h}$. Let us now compute the norm of $L=L_h$.

By definition, $\|L\|=\sup_{_{\|x\|_{_{L^2}}=1}} \|Lx\|_{_{L^{\infty}}}$. Now, for each $x\in L^2(\mathbb{R})$ we have that
\[
|Lx(t)|=\left | \int_{-\infty}^{\; \infty}x(s)\,h(t-s)\,ds\right |\leq \|x\|_{_{L^2}}\,\|h(t-\cdot)\|_{_{L^2}}=\|x\|_{_{L^2}}\,\|h\|_{_{L^2}},
\]
so that $\|L\|\leq \|h\|_{_{L^2}}$.
On the other hand, taking $x(t)=\overline{h(-t)}$ we have that
\[
|Lx(0)|=\left | \int_{-\infty}^{\; \infty}\overline{h(-s)}\,h(-s)\,ds\right |=\int_{-\infty}^{\; \infty}|h(s)|^2\,ds=\|h\|_{_{L^2}}^2=\|x\|_{_{L^2}}\,\|h\|_{_{L^2}},
\]
and therefore
\[
\|L\|\, \geq \, \big{\|}L\left(\frac{x}{\; \, \, \|x\|_{_{L^2}}}\right)\big{\|}_{_{L^{\infty}}}\, \geq \, \frac{|Lx(0)|}{\; \,\, \|x\|_{_{L^2}}}=\|h\|_{_{L^2}}.
\]
Hence $\|L\|=\|h\|_{_{L^2}}$. This proves the theorem. {\hfill $\Box$}

The following theorem  can be proved with the very same technique.

\begin{theorem}\label{two}
The map $\varphi:\ell^2(\mathbb{Z})\longrightarrow \mathbf{F_d}$ given by $\varphi(h)=L_h$, where $L_h(x)=x\ast h$, is an isometry. In particular, $\mathbf{F_d}$, with its usual norm (inherited from $\mathbf{B_d}$) is a Hilbert space with inner product given by $\langle L_h,L_g \rangle = \langle h,g \rangle$.
\end{theorem}

These two results allow us a perfect and natural identification between filters sending finite energy signals to bounded signals and the functions defining them as convolution operators. In particular, it is possible to compute the distance between two given filters (in the sense of their norms) as the distance between their unit impulse responses. Moreover, the computations can be made both in the time domain and in the frequency domain since the Fourier transform preserves norms (except for a multiplicative constant).

As we have already said, in this paper we will pay special attention to physi\-cally realizable filters and causal filters. Thanks to the previous theorems, these filters are easily characterized in the time domain. Concretely, the filter $L\in\mathbf{F_a}$ given by $L=L_{h}$ is causal (a physically realizable filter with time delay $T>0$) if and only if $h_{|_{(-\infty,0)}}= 0$ ($h_{|_{(-\infty,-T)}}=0$, respectively). Analogously, the filter $L\in\mathbf{F_d}$ given by $L=L_{h}$ is causal (a physically realizable filter with time delay $N>0$) if and only if we have that $h[n]=0$ for all $n\in\mathbb{Z}$ with $n<0$ ($n<-N$, respectively). We denote these classes of filters as $\mathbf{CF_a}$, $\mathbf{D}_T$,  $\mathbf{CF_d}$, and $\mathbf{DF}_N$, respectively.

\section{The analog case}
In this section we will use the notation $L_{[a,b]}$ for the analog filter with transfer function $\chi_{_{[a,b]}}$. We would like to estimate the distance $d(L_{[a,b]},\mathbf{CF_a})$ and the angle determined by $L_{[a,b]}$ and $\mathbf{CF_a}$.

We need first to recall some concepts and notations. The Fourier transform is the operator $\mathcal{F}:L^2(\mathbb{R})\longrightarrow L^2(\mathbb{R})$ given by $$\mathcal{F}(x)(w)=\frac{1}{\sqrt{2\pi}}\int_{-\infty}^{\; \infty}x(t)\,e^{-iwt}\,dt.$$ The inverse Fourier transform is $$\mathcal{F}^{-1}(x)(t)=\frac{1}{\sqrt{2\pi}}\int_{-\infty}^{\; \infty}x(w)\,e^{iwt}\,dw.$$ These operators are isometries.

Let us start by computing the function $h=\mathcal{F}^{-1}(\chi_{_{[a,b]}})$ (obviously, $L_{[a,b]}(x)=x\ast h$).
\begin{eqnarray*}
h(t) &=& \frac{1}{\sqrt{2\pi}}\int_{a}^{b}e^{iwt}\,dw=\left.\frac{e^{iwt}}{it\sqrt{2\pi}}\,\right]_{w=a}^{w=b}\\[2mm]
&=& \frac{1}{t\sqrt{2\pi}}\big[\, \sin(bt)-\sin(at)-i(\cos(bt)-\cos(at))\, \big]
\end{eqnarray*}

Clearly, the function $\tau(t)=|h(t)|^2=\frac{1}{\pi t^2}\left[1-\cos((b-a)t)\right]$ satisfies $\tau(t)=\tau(-t)$. Let us define
\[
h_{-}(t)=h(t)\, \chi_{_{(-\infty,0)}}(t) \quad \text{and} \quad h_+(t)=h(t)\, \chi_{_{[0,\infty)}}(t).
\]
Then $L_{[a,b]}=L_{h}=L_{h_{-}}+L_{h_{+}}$, $L_{h_{+}}\in\mathbf{CF_a}$ \, and
\[
\langle L_{h_{-}},L_{h_{+}}\rangle = \langle h_{-},h_{+}\rangle = \int_{-\infty}^{\; \infty}h_{-}(t)\,\overline{h_{+}(t)}\,dt=0.
\]
Hence $L_{h_{+}}$ is the best approximation to the ideal filter $L_{[a,b]}$ by elements of $\mathbf{CF_a}$ and the corresponding distance is given by:
\begin{eqnarray*}
d^2&=&d(L_{[a,b]},\mathbf{CF_a})^2=\|L_{h_{-}}\|^2=\|h_{-}\|^2  \\[3mm] &=& \int_{-\infty}^{\, 0}|h(t)|^2\,dt=\frac{1}{2}\int_{-\infty}^{\; \infty}|h(t)|^2\,dt=\frac{\|h\|^2}{2}=\frac{b-a}{2}.
\end{eqnarray*}
Moreover, now is obvious that the the angle determined by $L_{[a,b]}$ and $\mathbf{CF_a}$ is
\[
\theta =\arcsin\frac{d}{\|h\|}=\arcsin\frac{1}{\sqrt{2}}=\frac{\pi}{4}.
\]

Thus we have proved the following result:
\begin{theorem}For every frequency interval $[a,b]\subset \mathbb{R}$, the ideal analog filter $L_{[a,b]}$ forms an angle with the space of analog causal filters $\mathbf{CF_a}$ which is equal to ${\pi}/4$. Moreover (and in accordance with this phenomenon) the distance
$d=d(L_{[a,b]},\mathbf{CF_a})$ is $d=\sqrt{\frac{b-a}{2}}=\frac{1}{\sqrt{2}}\|L_{[a,b]}\|$.
\end{theorem}

Note that the angle does not depend on the frequencies interval $[a,b]$ and the distance depends just on its size. One may wonder why this is so. The natural answer is that the frequency response $H(\xi)$ of any ideal filter satisfies $H(\mathbb{R})\subseteq \mathbb{R}$, so that it satisfies the assumption of the following general theorem:

\begin{theorem}
Let us assume that $H(\mathbb{R})\subseteq \mathbb{R}$, where $H$ is the transfer function of the analog filter $L:L^2(\mathbb{R})\longrightarrow L^{\infty}(\mathbb{R})$. Then $d(L,\mathbf{CF_a})=\frac{1}{\sqrt{2}}\|L\|$ and the angle determined by $L$ and $\mathbf{CF_a}$ is ${\pi}/4$.
\end{theorem}

\noindent \textbf{Proof.} Set $h=\mathcal{F}^{-1}(H)$. Then
\[
A(t):=\mathbf{Re}(h(t))=\frac{1}{\sqrt{2\pi}}\int_{-\infty}^{\; \infty}H(\xi)\cos(t\xi)\,d\xi
\]
and
\[
B(t):=\mathbf{Im}(h(t))=\frac{1}{\sqrt{2\pi}}\int_{-\infty}^{\; \infty}H(\xi)\sin(t\xi)\,d\xi
\]

\noindent satisfy $A(-t)=A(t)$ and $B(-t)=-B(t)$, so that $\varphi(t)=|h(t)|^2$ is an even function and we can apply to $L$ the same decomposition technique we used for the ideal filter $L_{[a,b]}$. {\hfill $\Box$}
\medskip

Let us now compute the distance of (and the angle determined by) the ideal filter $L_{[a,b]}$ to the space $\mathbf{D}_T$ of physically realizable filters after a delay time $T>0$. These numbers can be computed by the same procedure, the main idea being to use the decomposition $h=g_{-}+g_{+}$ where
\[
g_{-}(t)=h(t)\, \chi_{_{(-\infty,-T)}}(t) \quad \text{and} \quad g_+(t)=h(t)\, \chi_{_{[-T,\infty)}}(t).
\]
The distance is now given by $$d(T)=\|g_{-}\|=\big\{\int_{-\infty}^{-T}\frac{1}{\pi t^2}\left[1-\cos((b-a)t)\right]dt\,\big\}^{\frac{1}{2}}.$$ Of course, we have that $$\|h\|^2=b-a=2\,d(T)^2+\int_{-T}^{\, T}\frac{1}{\pi t^2}\left[1-\cos((b-a)t)\right]dt.$$ Hence
\[
d(T)=\Big(\,\frac{b-a}{2}-\frac{1}{2}\int_{-T}^{\, T}\frac{1}{\pi t^2}\left[1-\cos((b-a)t)\right]dt\,\Big)^{\frac{1}{2}}.
\]
Now it is not difficult to prove the following result:
\begin{theorem}Let $d(T)$ and $\theta(T)$ denote the distance and the angle determined by $L_{[a,b]}$ and $\mathbf{D}_T$, respectively. Then
\vspace*{-2mm}
\begin{enumerate}
\item[{\rm (a)}] For every fixed interval $[a,b]$ and every $T>0$, $\theta(T)\in [0,\frac{\pi}{4}]$,
\vspace*{2mm}
\item[{\rm(b)}] For any fixed interval $[a,b]$, $\lim_{_{T\to \infty}}d(T)=\lim_{_{T\to \infty}}\theta(T)=0$,
\vspace*{2mm}
\item[{\rm(c)}] For any fixed $T>0$, there exists a constant $c(T)\in [0,\frac{2\pi}{T}]$ such that  $\lim_{_{(b-a)\to \infty}}d(T)=c(T)$. Moreover, $\lim_{_{(b-a)\to \infty}}\theta(T)=0$,
\vspace*{2mm}
\item[{\rm (d)}] For any fixed $T>0$, $\lim_{_{(b-a)\to 0}}d(T)=0$ and $\lim_{_{(b-a)\to 0}}\theta(T)=\frac{\pi}{4}$.
\end{enumerate}
\end{theorem}
Note that if we allow any positive delay time $T>0$, the angle determined by the ideal filters with a large bandpass interval is almost zero, what is in clear contrast with the result for causal filters (i.e., for $T=0$).
\medskip

\noindent \textbf{A philosophical remark} One may wonder why the angles we have got belong to the interval $[0,\pi/4]$. For example, why  ideal filters are never orthogonal to the causal ones? The natural answer is: The filter $L_h$ is orthogonal to $\mathbf{CF}_a$ if and only if it is a filter without memory! (Reason: For the filter $L_{h}$, being without memory means that $h_{_{|(0,\infty)}}\equiv 0$, which is equivalent to claim that $h\perp L^2(0,\infty)$\,). Thus, in a ``philosophical" sense, the angle between $L$ and $\mathbf{CF}_a$ is zero if and only if the filter is causal, the angle is $\pi/2$ if and only if the filter has no memory and the angle $\pi/4$ may be related to the fact that the filter uses all values (past and future) of the input signal. Of course, this is the case for the ideal filters and, in particular, we have shown that if we want to isolate exactly one frequency component of the signal, the associated filter determines an angle of $\pi/4$ with $\mathbf{D}_T$ for every $T\geq 0$.

\section{The digital case}
Let us consider a digital filter $L:\ell^2(\mathbb{Z})\longrightarrow\ell^{\infty}(\mathbb{Z})$. Then $Lx=x\ast h$ for a certain sequence $h=\{h(k)\}_{_{k=-\infty}}^{\infty}\in \ell^2(\mathbb{Z})$. In this case, the discrete Fourier transform $\mathcal{F}$ is given by
\[
\mathcal{F}_d\big(\{x[n]\}_{_{n=-\infty}}^{\infty}\big)(w)=\sum_{n=-\infty}^{\infty}x[n]\,e^{-iwn}=X(w),
\]
and satisfies
\begin{equation}\label{transformisometry}
\big{\|}\mathcal{F}_d\big(\{x[n]\}_{_{n=-\infty}}^{\infty}\big)(w)\big{\|}_{_{L^2(0,2\pi)}}^2 = 2\pi\, \big{\|}\{x[n]\}_{n=-\infty}^{\infty}\big{\|}_{_{\ell^2}}^2.
\end{equation}
Moreover,
\[
\mathcal{F}_d(x\ast h)(w)=\mathcal{F}_d(x)(w)\,\mathcal{F}_d(h)(w)=X(w)H(w),
\]
so that we can think about the ideal filter with bandpass given by $[a,b]$, $0<a<b<2\pi$, as the digital filter $L_{[a,b]}x=x\ast h_{[a,b]}$ with transfer function
$H_{[a,b]}(w)=\mathcal{F}_d(h_{[a,b]})(w)=\chi_{_{[a,b]}}(w)$. Unfortunately, this filter is not physically realizable.

Digital causal filters $L:\ell^2(\mathbb{Z})\longrightarrow\ell^{\infty}(\mathbb{Z})$  are characterized in the frequency domain
precisely as those with transfer function a $2\pi$-periodic function $H(w)\in \mathbb{CF}_d:=\mathbf{span}\{e^{-i w k}:k\in \mathbb{N}\}\subseteq L^2(0,2\pi)$ (here $\mathbf{span}$ denotes the closure of the subspace generated by $\{e^{-i w k}:k\in \mathbb{N}\}$ in $L^2(0,2\pi)$\,).
It follows from  Theorem \ref{two} and equation (\ref{transformisometry}) that
\[
d(L_{[a,b]},\mathbf{CF}_d)=\frac{1}{\sqrt{2\pi}}\, d(\chi_{[a,b]},\mathbb{CF}_d).
\]

Obviously, $\chi_{_{[a,b]}}\in L^2(0,2\pi)$ so that, denoting by $\{c_k\}$ its Fourier coefficients, $\chi_{_{[a,b]}}(w)=\sum_{k=-\infty}^{\infty}c_k \,e^{i k w}$. Moreover, $H(w)=\sum_{k=0}^{\infty}c_{-k}\, e^{-i k w}$ is the best approximation of $\chi_{_{[a,b]}}(w)$ in the $L^2$-norm by elements of $\mathbb{CF}_d$. Therefore,
\[
d^2:= d(\chi_{_{[a,b]}},\mathbb{F}_d)^2=\|\,\chi_{[a,b]}-H(w)\,\|^2=2\pi\sum_{k=1}^{\infty}|c_k|^2.
\]
Now, $$b-a=\|\chi_{_{[a,b]}}\|^2=2\pi|c_0|^2+2d^2$$ since $|c_k|=|c_{-k}|$ because $\chi_{_{[a,b]}}(\mathbb{R})\subseteq \mathbb{R}$. Moreover, $|c_0|^2=\frac{1}{4\pi^2}(b-a)^2$, and hence
\begin{equation}\label{distance_digital}
d(L_{[a,b]},\mathbf{CF}_d)=\frac{d}{\sqrt{2\pi}}=\sqrt{\frac{b-a}{2\pi}}\sqrt{\frac{1}{2}-\frac{b-a}{4\pi}},
\end{equation}

\begin{equation}\label{angle_digital}
\theta =\arcsin\sqrt{\frac{1}{2}-\frac{b-a}{4\pi}}.
\end{equation}
Thus, we have proved:
\begin{theorem}The angle $\theta$ determined by the ideal filter $L_{[a,b]}$ and the space of digital causal filters $\mathbf{CF}_d$, and the distance
$d(L_{[a,b]},\mathbf{CF}_d)$ are given by formulas $(\ref{angle_digital})$ and $(\ref{distance_digital})$, respectively. In particular, they are both functions of the length of the bandpass interval $[a,b]$. Moreover they satisfy the following properties:
\vspace*{-2mm}
\begin{enumerate}
\item[{\rm (a)}] $0<\theta<\pi/4$, \; $\lim_{_{(b-a)\to 0}}\theta=\pi/4$, \; $\lim_{_{(b-a)\to 2\pi}}\theta =0$,
\vspace*{2mm}
\item[{\rm (b)}]  $\lim_{_{(b-a)\to 0}}d(L_{[a,b]},\mathbf{CF}_d)=\lim_{_{(b-a)\to 2\pi}}d(L_{[a,b]},\mathbf{CF}_d) =0$.
\end{enumerate}
\end{theorem}

Note that the distance  and the angle do not depend on the location of the frequencies interval $[a,b]$ but just on its size. Moreover, the angle is never greater than $\pi/4$.
\medskip

\noindent \textbf{Remark} If $L\in \mathbf{F_d}$ has a transfer function $H(w)$ with the property that $H(\mathbb{R})\subseteq \mathbb{R}$, then $\|H\|^2=2\pi|c_0(H)|^2+2d(H,\mathbb{F}_d)^2$. Hence, if $|c_0(H)|=\frac{C}{\sqrt{2\pi}}\|H\|$ for a certain constant $C$, the angle determined by  $L$ and the space $\mathbf{CF}_d$ is
\[
\theta=\arcsin\sqrt{\frac{1-C^2}{2}}.
\]
In particular, if $H$ is an odd function then $\theta=\pi/4$.
\medskip

Let us now compute the distance of (and the angle formed by) the ideal filter $L_{[a,b]}$ to the space $\mathbf{DF}_N$ of physically realizable (after a delay time $N>0$) digital filters. First of all, we define
\[
 \mathbb{DF}_N=\mathcal{F}\big(\{h\in \ell^2(\mathbb{Z}):L_h\in \mathbf{DF}_N\}\big).
 \]
Clearly, $d(\chi_{_{[a,b]}},\mathbb{DF}_N)^2=\sum_{k=N+1}^{\infty}|c_k(\chi_{_{[a,b]}})|^2=:d(N)^2$. Now,
\begin{eqnarray*}
\|\chi_{_{[a,b]}}\|^2 &=& b-a=2\pi\sum_{k=-\infty}^{\infty}|c_k(\chi_{_{[a,b]}})|^2=2\,d(N)^2+2\pi\sum_{k=-N}^N|c_k(\chi_{_{[a,b]}})|^2\\
&=& 2\,d(N)^2+2\pi\,|c_0(\chi_{_{[a,b]}})|^2+2\sum_{k=1}^N2\pi\,|c_k(\chi_{_{[a,b]}})|^2.
\end{eqnarray*}
We know that $2\pi\,|c_0(\chi_{[a,b]})|^2=(b-a)^2/2\pi$. Moreover, an easy computation shows that, for $k\neq 0$,  $2\pi\,|c_k(\chi_{_{[a,b]}})|^2=\big(1-\cos(k(b-a))\big)\,/k^2\pi$, so that
\begin{equation}\label{distance}
d(L_{[a,b]},\mathbf{DF}_N)=\frac{d(N)}{\sqrt{2\pi}}=\sqrt{\frac{b-a}{2\pi}}\left \{ \frac{1}{2}-\frac{b-a}{4\pi}-\sum_{k=1}^N\frac{1-\cos(k(b-a))}{k^2\pi(b-a)}\right\}^{\frac{1}{2}},
\end{equation}
and the angle determined by $L_{[a,b]}$ and $\mathbf{DF}_N$ is
\begin{equation}\label{angle}
\theta(N) =\arcsin\left \{ \frac{1}{2}-\frac{b-a}{4\pi}-\sum_{k=1}^N\frac{1-\cos(k(b-a))}{k^2\pi(b-a)}\right\}^{\frac{1}{2}}.
\end{equation}
Thus the following result holds true:
\begin{theorem}The angle $\theta(N)$ determined by the ideal filter $L_{[a,b]}$ and the space  $\mathbf{DF}_N$,  and the distance
$d(L_{[a,b]},\mathbf{DF}_N)$ are given by formulas $(\ref{angle})$ and $(\ref{distance})$, respectively. In particular, they are both functions of the length of the bandpass interval $[a,b]$. Moreover they satisfy the following properties:
\vspace*{-2mm}
\begin{enumerate}
\item[{\rm (a)}] $0<\theta(N)<\pi/4$, \; $\lim_{_{(b-a)\to 0}}\theta(N)=\pi/4$, \; $\lim_{_{(b-a)\to 2\pi}}\theta(N) =0$,
\vspace*{2mm}
\item[{\rm (b)}]  $\lim_{_{(b-a)\to 0}}d(L_{[a,b]},\mathbf{DF}_N)=\lim_{_{(b-a)\to 2\pi}}d(L_{[a,b]},\mathbf{DF}_N) =0$,
\vspace*{2mm}
\item[{\rm (c)}] Given the bandpass interval $[a,b]\subset (0,2\pi)$, the sequence $\{\theta(N)\}_{_{N=0}}^{\infty}$ is strictly decreasing and $\lim_{_{N\to \infty}}\theta(N)=0$.
\end{enumerate}
\end{theorem}

\section*{Acknowledgements}
 The first author has been partially supported by Research Project P07-TIC-02713 and Research Group FQM-332. The second author has been partially
 supported by MCYT-FEDER Grant BFM2001-2871-C04-01.

\medskip

Received November 30, 2008; revised December 2008.

\medskip

\end{document}